\documentclass[preprint,12pt]{elsarticle}

\input epsf
\usepackage{amssymb}

\newcounter{bla}

\def\CHF#1#2#3{
{}_1F_1\left(
\begin{array}{c}
\begin{array}{cc} \hskip-10pt#1 \end{array}\\
\begin{array}{c}  \hskip-10pt#2 \end{array}
\end{array}
\hskip-8pt;\,#3
\right)}

\def\Ai{{{\rm Ai}}}

\def\arcsinh{{\rm arcsinh}}
\def\arccosh{{\rm arccosh}}

\def\sk{\sum_{k=0}^\infty\,}

\def\binomial#1#2{
\renewcommand{\arraystretch}{1.0}
\left(
\begin{array}{c} 
\hskip-5pt#1\\
\hskip-5pt#2
\end{array}
\hskip-5pt\right)}
\def\eqref#1{(\ref{#1})}

\def\tfrac#1#2{{{\lower.6ex
\hbox{$\scriptstyle#1$}}\over
{\raise.7ex
\hbox{$\scriptstyle#2$}}}}

\def\phase{{\rm ph}}
\def\bigO{{\cal O}}

\def\arccosh{{\rm arccosh}}

\def\dsp#1{\displaystyle#1}

\def\dsp{\displaystyle}
\def\Frac#1#2{\frac{\displaystyle{#1}}{\displaystyle{#2}}}

\def\bigO{{\cal O}}

\begin{document}

\begin{frontmatter}



\title{Efficient computation of Laguerre polynomials}


\author[1,3]{Amparo Gil}
\author[2,3]{Javier Segura}
\author[4]{Nico M. Temme}
\address[1]{Depto. de Matem\'atica Aplicada y Ciencias de la Comput. Universidad de Cantabria. 39005-Santander, Spain. e-mail: amparo.gil@unican.es}
\address[2]{Depto. de Matem\'aticas, Estad\'{\i}stica y Comput. Universidad de Cantabria. 39005-Santander, Spain. e-mail: javier.segura@unican.es  }
\address[3]{Department of Mathematics and Statistics. San Diego State University. 5500 Campanile Drive San Diego, CA, USA.}
\address[4]{ IAA, 1825 BD 18, Alkmaar, The Netherlands\footnote{Former address: CWI, 1098 XG Amsterdam, The Netherlands}. e-mail: nico.temme@cwi.nl  }

\begin{abstract}
  An efficient algorithm and a Fortran 90 module ({\bf LaguerrePol}) for computing Laguerre polynomials $L^{(\alpha)}_n(z)$
  are presented. The standard three-term recurrence relation satisfied by the polynomials and different types of asymptotic expansions
  valid for $n$ large and $\alpha$ small, are used depending on the parameter region. 

Based on tests of contiguous relations in the parameter $\alpha$ and the degree $n$ satisfied by the polynomials, we claim that
a  relative accuracy close or better than $10^{-12}$ can be obtained using the module
{\bf LaguerrePol} for computing the functions $L^{(\alpha)}_n(z)$ in the parameter range $z \ge 0$, $-1 < \alpha \le 5$, $n \ge 0$. 
\end{abstract}
\end{frontmatter}

\section{Introduction}

  As is well known, Laguerre polynomials $L^{(\alpha)}_n(z)$ are involved in a vast number of
applications in physics (quantum mechanics, plasma physics, etc) and engineering; for example
see \cite{Kuang:1997:MIG,Bao:2008:GLH}. Also, the evaluation of Laguerre
 polynomials is central in the computation of nodes and weights in  Gauss-Laguerre quadrature rules. 

In this paper, we present an algorithm for computing Laguerre polynomials based on the use of
the standard three-term recurrence relation satisfied
by the polynomials and three types of asymptotic expansions valid for $n$ large and small values of the parameter $\alpha$: two Bessel-type
expansions (used in the oscillatory region of the functions) and a uniform Airy-type expansion. This Airy expansion is specially suitable in the transition 
 between the oscillatory and monotonic regimes of the Laguerre polynomials, although its domain of applicability extends 
to a large part of both the oscillatory and the monotonic regions.  

The resulting algorithm, implemented in the Fortran 90 module {\bf LaguerrePol} is accurate and particularly efficient
for large values of the parameter $n$.

\section{Theoretical background}

Laguerre polynomials  $L^{(\alpha)}_n(z)$ are solutions of the differential equation

\begin{equation}
\label{odelag}
zy''+(\alpha + 1-z)y'+ny=0\,.
\end{equation}

The polynomials present an oscillatory and a monotonic regime, depending on the parameter values.
The oscillatory (monotonic) region of $L^{(\alpha)}_n(z)$  is found in the interval
$0<z/\nu<1$ ($z/\nu>1$), where 

\begin{equation}
\label{eq:nu}
\nu=4\left(n+\tfrac12(\alpha+1)\right). 
\end{equation}

Sharper lower ($z_l$)
and upper ($z_u$) bounds 
limiting the region where the 
Laguerre polynomial oscillates are given by \cite{Dimitrov:2010:BEZ}:

\begin{equation}\label{eq:lagairy1}
\begin{array}{lcl}
 z_l&=&\Frac{2n^2+n(\alpha-1)+2(\alpha+1)-2(n-1)\displaystyle{\sqrt{n^2+(n+2)(\alpha+1)}}}{(n+2)},\\
 z_u&=&\Frac{2n^2+n(\alpha-1)+2(\alpha+1)+2(n-1)\displaystyle\sqrt{n^2+(n+2)(\alpha+1)}}{(n+2)}.
\end{array}
\end{equation}

An example of the behaviour of the Laguerre polynomials  in the oscillatory region is shown in
Figure 1.  

\begin{figure}
\begin{center}
\epsfxsize=13cm \epsfbox{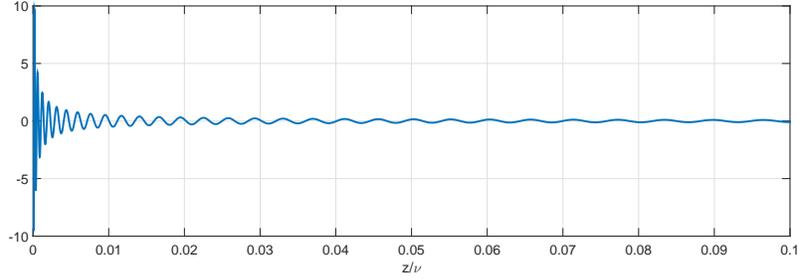}
\end{center}
\caption{The function ${ L}^{(4)}_{150}(z)$ is plotted as an example. The value of the parameter $\nu$ (defined
in (\ref{eq:nu})) is 610 in this particular case.
\label{Fig1}}
\end{figure}

Next, we are going to describe the theoretical expressions involved in the 
computation of Laguerre polynomials.

\subsection{Asymptotic expansions for $n$ large}

In the asymptotic expansions for large degree we assume the $\alpha$ is fixed, which means small with respect to $n$. If we need results for which $\alpha$  is not small enough, we can use recursion with respect to $\alpha$. That is,
\begin{equation}\label{eq:lagasy01}
xL_n^{(\alpha+1)}(x)=(\alpha+x)L_n^{(\alpha)}(x)-(\alpha+n)L_n^{(\alpha-1)}(x).
\end{equation}
This  follows from the corresponding $c$-recursion of the Kummer function $U(-n,c,x)=(-1)^n\, n!\, L_n^{(c-1)}(x)$.

 \subsubsection{An expansion in terms of Airy functions}\label{sec:LagAiryexp}

We start with the representation\footnote{We summarize results of \cite{Frenzen:1988:UAE};  see also \cite[\S{VII.5}]{Wong:2001:AAI}.}
\begin{equation}\label{eq:lagairy01}
L_n^{(\alpha)}(\nu x)=(-1)^n
\frac{ e^{\frac12\nu x}\chi(\zeta)}{2^\alpha\nu^{\frac13}}\left(\Ai\left(\nu^{2/3} \zeta\right)A(\zeta)
+\nu^{-\frac43}\Ai^{\prime}\left(\nu^{2/3}\zeta\right)
B(\zeta)\right)
\end{equation}
with  expansions
\begin{equation}\label{eq:lagairy02}
A(\zeta)\sim\sum_{j=0}^\infty\frac{\alpha_{2j}}{\nu^{2j}},\quad B(\zeta)\sim\sum_{j=0}^\infty\frac{\beta_{2j+1}}{\nu^{2j}},\quad n\to\infty,
\end{equation}
uniformly for bounded $\alpha$ and  $x\in(x_0,\infty)$, where $x_0\in(0,1)$, a fixed number. 

Here
\begin{equation}\label{eq:lagairy03}
\nu=4\kappa,\quad \kappa=n+\tfrac12(\alpha+1), \quad 
\chi(\zeta)=2^{\frac12}x^{-\frac14-\frac12\alpha}\left(\frac{\zeta}{x-1}\right)^{\frac14},
\end{equation}
and

\begin{equation}
 \label{eq:lagairy04} 
\begin{array}{ll}
\dsp{\tfrac23(-\zeta)^{\frac32}=\tfrac12\left(\arccos\sqrt{{x}}-\sqrt{{x-x^2}}\right)}  & \mbox{if $ 0<x \le 1$,}\\
\dsp{\tfrac23\zeta^{\frac32}=\tfrac12\left(\sqrt{{x^2-x}}-\arccosh\sqrt{{x}}\right)} & \mbox{if $ x\ge 1$.}
\end{array}
\end{equation}
We have the relation
\begin{equation}\label{eq:lagairy05}
\zeta^\frac12\frac{d\zeta}{dx}=\frac{\sqrt{x-1}}{2\sqrt{x}}.
\end{equation}
%

%
%

The first coefficients of the expansions in \eqref{eq:lagairy02} are
\begin{equation}\label{eq:lagairy07}
\alpha_0=1,\quad\ \beta_1=-\frac{1}{4b^3}\left(f_1-bf_2\right),
\end{equation}
where $b=\sqrt{\zeta}$ if $\zeta\ge0$ and $b=i\sqrt{-\zeta}$ when $\zeta\le0$, and
\begin{equation}\label{eq:lagairy08}
\begin{array}{ll}
\dsp{f_1=i\frac{\left(x+3\alpha(x-1)\right)x^2a_1^3-2}{3a_1^2x\sqrt{x(1-x)}},}\\[8pt] 
\dsp{f_2=\frac{-4-8x^2(x+3x\alpha-3\alpha)a_1^3+x^4(12x-3-4x^2+12\alpha^2(x-1)^2)a_1^6}{12x^3a_1^4(x-1)},}\\[8pt] 
\dsp{a_1=\left(\frac{4\zeta}{x^3(x-1)}\right)^{\frac14}}.
\end{array}
\end{equation}
More coefficients can be obtained by the method described in \cite[\S23.3]{Temme:2015:AMI}.

\subsubsection{A simple Bessel-type expansion}\label{sec:Lagsim}

For the Laguerre polynomials we consider two types of asymptotic expansions in terms of Bessel functions, 
one for small values of the variable $z$ of $L_n^{(\alpha)}(z)$ and one in which larger values are allowed.
Next we give some details of the asymptotic expansion valid for small values of the variable $z$.

We use the expansion of the Kummer function ${}_1F_1(a;c;x)$ for large negative values of $a$. First we mention
\begin{equation}\label{eq:lagsim01}
L_n^{(\alpha)}(z)=\binomial{n+\alpha}{n}\CHF{-n}{\alpha+1}{z}.
\end{equation}
Then, see  \cite[\S10.3.4]{Temme:2015:AMI},
\begin{equation}\label{eq:lagsim02}
\begin{array}{ll}
\dsp{\frac{1}{\Gamma(c)}\CHF{-a}{c}{z}\sim\left(\frac{z}{a}\right)^{\frac12(1-c)}\frac{\Gamma(1+a)e^{\frac12z}}{\Gamma(a+c)}\ \times}\\[8pt]
\quad\quad\dsp{\left(J_{c-1}\left(2\sqrt{az}\right)\sk \frac{a_k(z)}{(-a)^k}-\sqrt{\frac{z}{a}}J_{c}\left(2\sqrt{az}\right)\sk \frac{b_k(z)}{(-a)^k}\right).}
\end{array}
\end{equation}
This expansion of ${}_1F_1(-a;c;z)$ is valid for bounded values of $z$ and $c$, with $a\to\infty$ inside the sector  $-\pi+\delta\le\phase\,a\le\pi-\delta$. This gives for the Laguerre polynomial
\begin{equation}\label{eq:lagsim03}
\begin{array}{ll}
\dsp{L_n^{(\alpha)}(x)\sim\left(\frac{x}{n}\right)^{-\frac12\alpha} e^{\frac12x}\ \times}\\[8pt]
\quad\dsp{\left(J_{\alpha}\left(2\sqrt{nx}\right)\sk (-1)^k\frac{a_k(x)}{n^k}-\sqrt{\frac{x}{n}}J_{\alpha+1}\left(2\sqrt{nx}\right)\sk (-1)^k \frac{b_k(x)}{n^k}\right), \quad n\to\infty.}
\end{array}
\end{equation}

The coefficients $a_k(x)$ and $b_k(x)$ follow from the expansion of the function
\begin{equation}\label{eq:lagsim04}
f(z,s)=e^{xg(s)}\left(\frac{s}{1-e^{-s}}\right)^{\alpha+1},\quad g(s)=\frac{1}{s}-\frac{1}{e^s-1}-\frac12.
\end{equation}
The function $f$ is analytic in the strip $\vert\Im s\vert<2\pi$ and it can be expanded for  $\vert s\vert<2\pi$ into
\begin{equation}\label{eq:lagsim05}
f(x,s)=\sk c_k(x) s^k.
\end{equation}
The coefficients $c_k(x)$  are combinations of Bernoulli numbers and Bernoulli polynomials, the first ones being (with $c=\alpha+1$)
\begin{equation}\label{eq:lagsim06}
\begin{array}{@{}r@{\;}c@{\;}l@{}}
c_0(x)&=&1,\quad c_1(x)=\frac{1}{12}\left(6c-x\right),\quad \\[8pt]
c_2(x)&=&\frac{1}{288}\left(-12c+36c^2-12xc+x^2\right),\\[8pt]
c_3(x)&=&\frac{1}{51840}\left(-5x^3 + 90x^2c +(-540c^2 + 
 180c+72)x +1080c^2(c-1)\right).
\end{array}
\end{equation}
The coefficients $a_k(x)$ and $b_k(x)$ are in terms of the $c_k(x)$ given by
\begin{equation}\label{eq:lagsim07}
\begin{array}{@{}r@{\;}c@{\;}l@{}}
a_k(x) & = & \dsp{\sum_{m=0}^k \binomial{k}{m}(m+1-c)_{k-m}x^m c_{k+m}(x),}\\[8pt]
b_k(x) & = & \dsp{\sum_{m=0}^k \binomial{k}{m}(m+2-c)_{k-m}x^m c_{k+m+1}(x),}
\end{array}
\end{equation}
$k=0,1,2,\ldots$, and the first relations are
\begin{equation}\label{eq:lagsim08}
\begin{array}{ll}
a_0(x)= c_0(x)=1,\quad b_0(x)= c_1(x),\\[8pt]
a_1(x)= (1-c)c_1(x)+xc_2(x),\quad b_1(x)= (2-c)c_2(x)+xc_3(x),\\[8pt]
a_2(x)= (c^2-3c+2)c_2(x)+(4x-2xc)c_3(x)+x^2c_4(x),\\[8pt]
b_2(x)=  (c^2-5c+6)c_3(x)+(6x-2xc)c_4(x)+x^2c_5(x),
\end{array}
\end{equation}
again with $c=\alpha+1$.

\subsubsection{A not so simple expansion in terms of Bessel functions}\label{sec:LagBesexp}

In this case we use the representation\footnote{We summarize the results of \cite{Frenzen:1988:UAE}; see also \cite[\S{VII.7}]{Wong:2001:AAI}.} 
\begin{equation}\label{eq:lagbess01}
L_n^{(\alpha)}(2\nu x)= \frac{e^{\nu x}\chi(\zeta)}{2^\alpha\zeta^{\frac12\alpha}}\left(J_\alpha\bigl(2\nu \sqrt{\zeta}\bigr)
A(\zeta)-
\frac{1}{\nu\sqrt{\zeta}}J_{\alpha+1}\bigl(2\nu  \sqrt{\zeta}\bigr)B(\zeta)\right),
\end{equation}
with expansions
\begin{equation}\label{eq:lagbess02}
A(\zeta)\sim\sum_{j=0}^\infty\frac{A_j(\zeta)}{\nu^{2j}},\quad 
B(\zeta)\sim\sum_{j=0}^\infty\frac{B_j(\zeta)}{\nu^{2j}},\quad \nu\to\infty.
\end{equation}
Here,
\begin{equation}\label{eq:lagbess03}
\nu=2n+\alpha+1,\quad \chi(\zeta)=(1-x)^{-\frac14}\left(\frac{\zeta}{x}\right)^{\frac12\alpha+\frac14},\quad x<1,
\end{equation}
with $\zeta$ given by

\begin{equation}\label{eq:lagbess04}
\begin{array}{ll}
\dsp{\sqrt{-\zeta}=\tfrac12\left(\sqrt{{x^2-x}}+\arcsinh\sqrt{{-x}}\right)}, & \quad\mbox{if \quad$x\le0$,}\\[8pt]
\dsp{\sqrt{\zeta}=\tfrac12\left(\sqrt{{x-x^2}}+\arcsin\sqrt{{x}}\right)}, & \quad \mbox{if \quad$ 0\le x<1$.}
\end{array}
\end{equation}

 We have the relation
\begin{equation}\label{eq:lagbess05}
\frac{1}{\zeta^{\frac12}}\frac{d\zeta}{dx}=\sqrt{\frac{1-x}{x}},\quad x <1.
\end{equation}

The first coefficients are
\begin{equation}\label{eq:lagbess06}
\begin{array}{@{}r@{\;}c@{\;}l@{}}
A_0(\zeta)&=&1, \\[8pt] 
B_0(\zeta)&=&\dsp{\frac{1-4\alpha^2}{16}+\frac{\sqrt{\zeta}}{8\sqrt{\xi}}\left(\frac{4\alpha^2-1}{2}+\xi+\frac{5}{6}\xi^2\right),\quad
\xi=\frac{x}{1-x}.}
\end{array}
\end{equation}

We give a few details about the coefficients $A_j(\zeta)$ and $B_j(\zeta)$ of the expansions in \eqref{eq:lagbess02}. The first ones are given in \eqref{eq:lagbess06}. 

First we need coefficients $c_k^\pm$ of the expansions
\begin{equation}\label{eq:lagbess11}
s=s_+ + \sum_{k=1}^\infty c_k^+(u-ib)^k, \quad
s=s_- + \sum_{k=1}^\infty c_k^-(u+ib)^k,
\end{equation}
where, for $0\le x<1$, $b$, $s_\pm$ and the relation between $s$ and $u$ are defined by
\begin{equation}\label{eq:lagbess12}
b=\sqrt{\zeta},\quad s_\pm=\pm i\arcsin\sqrt{x},\quad s-x\coth s=u-b^2/s,
\end{equation}
with $\zeta$ defined in \eqref{eq:lagbess04}, and $s_\pm$ being the saddle points of the $s$-function and $\pm ib$ of the $u$-function.

Because $s(u)$ is an odd function of $u$, we have
\begin{equation}\label{eq:lagbess13}
c_k^-=(-1)^{k+1}c_k^+,\quad k=1,2,3,\ldots.
\end{equation}
In the following we write $c_k^+=c_k$. The first coefficients are
\begin{equation}\label{eq:lagbess14}
\begin{array}{@{}r@{\;}c@{\;}l@{}}
c_1&=&\dsp{\left(\frac{x}{\zeta(1-x)}\right)^{\frac14},}\\[8pt]
c_2&=&\dsp{-\frac{i{c_1}^2\left(3-3c_1-2x+3xc_1\right)}{6\sqrt{x(1-x)}},} \\[8pt]
c_3&=& \dsp{-\frac{c_1^3\bigl((27c_1^2-24c_1+8)x^2-(54c_1^2-60c_1+12)x+27c_1^2-36c_1+9\bigr)}{72x(1-x)}.}
\end{array}
\end{equation}

Next we consider the function with expansion
\begin{equation}\label{eq:lagbess15}
h(u)=\frac{1}{\chi(\zeta)}\left(\frac{\sinh s}{u}\right)^{-\alpha-1}\frac{ds}{du},\quad h(u)=\sum_{k=0}^\infty d_k(u-ib)^k,
\end{equation}
where $\chi(\zeta)$ is defined in \eqref{eq:lagbess03}. Because $h(u)$ is even, an expansion at the other saddle point $-ib$ has coefficients $(-1)^kd_k$; also, $d_0=h(ib)=1$. This makes $A_0(\zeta)=1$, see \eqref{eq:lagbess06}. The coefficient $d_1$ is given by
\begin{equation}\label{eq:lagbess16}
d_1=\frac{2ibc_2x+\gamma c_1x+i\gamma c_1^2b(1-x)+\left(2ibc_2-\gamma c_1^2b-i\gamma c_1\right)\sqrt{x(1-x)}}{c_1\sqrt{\zeta x}},
\end{equation}
where $\gamma=\alpha+1$.

The coefficients in \eqref{eq:lagbess02} follow from the following recursive scheme
\begin{equation}\label{eq:lagbess17}
h_k(u)=\alpha_k+\beta_k/u+\left(1+b^2/u^2\right)g_k(u), \quad
h_{k+1}(u)=g_k^\prime(u)-\frac{\alpha+1}{u}g_k(u),
\end{equation}
$k=0,1,2,\ldots$, with starting value $h_0(u)=h(u)$. The coefficients $\alpha_k$ and $\beta_k$ follow from substituting $u=\pm ib$. Because $h(u)$ is even, $\beta_0=0$, and $h_1(u)$ is odd, giving $\alpha_1=0$, and so on. Then, the coefficients in the expansions in \eqref{eq:lagbess02} are given by
\begin{equation}\label{eq:lagbess18}
A_k(\zeta)=\alpha_{2k},\quad B_k(\zeta)=\beta_{2k+1},\quad k=0,1,2,\ldots.
 \end{equation}
This gives, again with $\gamma=\alpha+1$,
\begin{equation}\label{eq:lagbess19}
\begin{array}{@{}r@{\;}c@{\;}l@{}}
A_0(\zeta)&=& 1, \quad
B_0(\zeta)= -\Frac{b}{4}\Bigl(3id_1+2i\gamma d_1+2d_2b\Bigr), \\[8pt]
A_1(\zeta)&=& -\Frac{1}{32b}\Bigl(48\gamma bd_2-46d_2b+24d_4b^3-8\gamma ^2bd_212i)\gamma ^2d_1\ +\\[8pt]
&&24ib^2\gamma d_3-24i\gamma d_1-60ib^2d_3+9id_1\Bigr),\\[8pt]
B_1(\zeta)&=& \Frac{1}{128b}\Bigl(-572ib^2\gamma d_3+36i\gamma ^2d_1-840ib^4d_5+480ib^2d_3\ +\\[8pt]
&&240ib^4\gamma d_5-46i\gamma d_1+192i\gamma ^2b^2d_3-8i\gamma ^3d_1-96b^3\gamma ^2d_4\ +\\[8pt]
&&15id_1-16b\gamma ^3d_2-16ib^2\gamma ^3d_3+72\gamma ^2bd_2+30d_2b-1032d_4b^3\ +\\[8pt]
&&240d_6b^5-92\gamma bd_2+672\gamma b^3d_4\Bigr).
\end{array}
 \end{equation}

For small values of $x$ we need expansions. We can expand in terms of $\zeta$ or $x$. For example, we can write
\begin{equation}\label{eq:lagbess21}
A_k(\zeta)=\sum_{j=0}^\infty A_j^{(k)} x^k, \quad B_k(\zeta)=\sum_{j=0}^\infty B_j^{(k)} x^k.
 \end{equation}
The first coefficients are $A_0^{(0)}=1$, $A_j^{(0)}=0$ $(j\ge1)$, $B_0^{(k)}=0$  $(k\ge0)$, and
\begin{equation}\label{eq:lagbess22}
\begin{array}{@{}r@{\;}c@{\;}l@{}}
A_0^{(1)}&=& -\frac{1}{6} \alpha(\alpha^2-1), \\[8pt]
A_1^{(1)}&=&-\frac{1}{360} (\alpha-1)(5\alpha-7)(\alpha+3)(\alpha+2),\\[8pt]
A_0^{(2)}&=&  \frac{1}{360}\alpha(\alpha^2-1)(\alpha-2)(\alpha-3)(5\alpha+7),\\[8pt]
A_1^{(2)}&=&\frac{1}{4536}(\alpha-1)(\alpha^2-4)(\alpha-3)(\alpha+3)(7\alpha^2-31),\\[8pt]
B_1^{(0)}&=&-\frac{1}{6} (\alpha^2-1),\\[8pt]
B_1^{(1)}&=&  \frac{1}{360}(\alpha^2-1)(\alpha-2)(\alpha-3)(5\alpha+7),\\[8pt]
B_1^{(2)}&=&  -\frac{1}{45360} (\alpha^2-1)(\alpha-2)(\alpha-3)(\alpha-4)(\alpha-5)(35\alpha^2+112\alpha+93).\\[8pt]
\end{array}
 \end{equation}

\subsubsection{Expansions for large  values of  $n$ and $\alpha$}\label{sec:LagBesexpla}
In \cite{Temme:1986:LLD} we have given expansions for large $n$ in which $\alpha=\bigO(n)$ is allowed; for a summary see 
\cite{Temme:1990:AEL}. These results can be obtained by using an integral representation, but they 
follow also from uniform expansions of
Whittaker functions obtained by using differential equations; see 
\cite{Dunster:1989:UAE}. These expansions include the $J$-Bessel function, and are valid in the parameter domain where order and argument of the Bessel function are equal, that is, in the turning point domain.  
Because no explicit forms of the coefficients in the expansions are available, we omit further details.

\subsubsection{An algorithm for computing the Bessel functions $J_{\nu}(z)$ }

The algorithm for computing the Bessel function $J_{\nu}(z)$ in the expansions (\ref{eq:lagsim03}) and (\ref{eq:lagbess01})
 is based in the following methods
of approximation:

\begin{description}

\item{\bf Power series.}

The power series given in Eq.(10.2.2) of \cite[\S10.19(ii)]{Olver:2010:Bessel}
is used for computing $J_{\nu}(z)$ when $z$ is small:

\[\mathop{J_{\nu}\/}\nolimits\!\left(z\right)=(\tfrac{1}{2}z)^{\nu}\sum_{k=0}^{%
\infty}(-1)^{k}\frac{(\tfrac{1}{4}z^{2})^{k}}{k!\mathop{\Gamma\/}\nolimits\!%
\left(\nu+k+1\right)}.\]

\item{\bf Debye's asymptotic expansions.}

Debye's asymptotic expansions  are also used in the algorithm.
The expressions are given in  Eq.(10.19.3) and Eq.(10.19.6) of
\cite[\S10.19(ii)]{Olver:2010:Bessel}: 

When $\nu <z$, we use

\[\mathop{J_{\nu}\/}\nolimits\!\left(\nu\mathop{\mathrm{sech}\/}\nolimits\alpha%
\right)\sim\frac{e^{\nu(\mathop{\tanh\/}\nolimits\alpha-\alpha)}}{(2\pi\nu%
\mathop{\tanh\/}\nolimits\alpha)^{\frac{1}{2}}}\sum_{k=0}^{\infty}\frac{U_{k}(%
\mathop{\coth\/}\nolimits\alpha)}{\nu^{k}},\]
and for $\nu >z$ 

\[\mathop{J_{\nu}\/}\nolimits\!\left(\nu\mathop{\sec\/}\nolimits\beta\right)\sim%
\left(\frac{2}{\pi\nu\mathop{\tan\/}\nolimits\beta}\right)^{\frac{1}{2}}\*%
\left(\mathop{\cos\/}\nolimits\xi\sum_{k=0}^{\infty}\frac{U_{2k}(i\mathop{\cot%
\/}\nolimits\beta)}{\nu^{2k}}-i\mathop{\sin\/}\nolimits\xi\sum_{k=0}^{\infty}%
\frac{U_{2k+1}(i\mathop{\cot\/}\nolimits\beta)}{\nu^{2k+1}}\right).\]

The coefficients $U_k(p)$ are polynomials in $p$ of degree $3k$ given by $U_0(p)=1$ and

\[U_{k+1}(p)=\tfrac{1}{2}p^{2}(1-p^{2})U_{k}^{\prime}(p)+\frac{1}{8}\int_{0}^{p}%
(1-5t^{2})U_{k}(t)dt.\]

\item{\bf Asymptotic expansions for large $z$.}

 For large values of the argument $z$, we use the Hankel's expansion given
in  \cite[\S10.17(i)]{Olver:2010:Bessel}:

\[\mathop{J_{\nu}\/}\nolimits\!\left(z\right)\sim\left(\frac{2}{\pi z}\right)^{%
\frac{1}{2}}\*\left(\mathop{\cos\/}\nolimits\omega\sum_{k=0}^{\infty}(-1)^{k}%
\frac{a_{2k}(\nu)}{z^{2k}}-\mathop{\sin\/}\nolimits\omega\sum_{k=0}^{\infty}(-%
1)^{k}\frac{a_{2k+1}(\nu)}{z^{2k+1}}\right),\]
where

\[\omega=z-\tfrac{1}{2}\nu\pi-\tfrac{1}{4}\pi.\]

The coefficients $a_{k}(\nu)$ are given by

\[a_{k}(\nu)=\frac{(4\nu^{2}-1^{2})(4\nu^{2}-3^{2})\cdots(4\nu^{2}-(2k-1)^{2})}{%
k!8^{k}}.\]

\item{\bf Airy-type expansions.}

An important ingredient in our algorithm for computing Bessel functions are 
Airy-type expansions.
We use the representation given in \cite[Chapter~8]{Gil:2007:NSF} 

$$
\dsp{J_\nu(\nu x)
=\quad \frac{\phi(\zeta)}{\nu^{1/3}}\left[\Ai(\nu^{2/3}\zeta)\,{A}_{\nu}(\zeta)
+\nu^{-4/3}\Ai'(\nu^{2/3}\zeta)\,{B}_{\nu}(\zeta)\right],}\\ 
\renewcommand{\arraystretch}{1.0}
$$

where 

$$
\phi(\zeta)=\left(\frac{4\zeta}{1-x^2}\right)^{\frac14},\quad \phi(0)=2^{\frac13}.
$$

The variable $\zeta$ is written in terms of the variable $x$ as

\renewcommand{\arraystretch}{1.5}
$$
\begin{array}{llll}
\dsp{\tfrac23\zeta^{3/2}}&=&\dsp{\ln\frac{1+\sqrt{{1-x^2}}}x-\sqrt{{1-x^2}},} &  \quad 0\le x\le1,\\
\dsp{\tfrac23(-\zeta)^{3/2}}&=&\dsp{\sqrt{{x^2-1}}-\arccos\frac{1}{x},} & \quad x\ge1.
\end{array}
\renewcommand{\arraystretch}{1.0}
$$

\item{\bf Three-term recurrence relation using Miller's algorithm.}

The standard three-term recurrence relation for the cylinder functions  

\[\mathop{{J}_{\nu-1}\/}\nolimits\!\left(z\right)+\mathop{{J}_{%
\nu+1}\/}\nolimits\!\left(z\right)=(2\nu/z)\mathop{{J}_{\nu}\/}%
\nolimits\!\left(z\right),\]

is computed backwards (starting from large values of $\nu$) using Miller's algorithm.

\end{description}

\subsection{Three-term recurrence relation}

The generalized Laguerre polynomials satisfy the following three-term recurrence relation

\begin{equation}
\label{ttrr}
L^{(\alpha)}_{n+1}(z)=\Frac{2n+\alpha+1-z}{n+1}L^{(\alpha)}_n(z)-\Frac{n+\alpha}{n+1}
L^{(\alpha)}_{n-1}(z)\,.
\end{equation}

This recurrence relation is not ill conditioned in both backward and forward directions.
 Therefore it can be used with starting values $L^{(\alpha)}_0(z)=1$
and $L^{(\alpha)}_1(z)=1+\alpha-z$, to compute 
the generalized Laguerre polynomials when $n$ is small/moderate. As $n$ increases,
it is more efficient, as we later discuss, to use the asymptotic expansions described in the previous sections.

\section{Overview of the software structure}

The Fortran 90 package includes the main module {\bf LaguerrePol}, which includes
as public routine the function {\bf laguerre}.

In the module {\bf LaguerrePol}, the auxiliary 
modules {\bf Someconstants}  (a module for the 
computation of the main constants used in 
  the different routines),  {\bf BesselJY} (for the computation of Bessel functions) and {\bf AiryFunction} 
(for the computation of Airy functions) are used.

\section{Description of the individual software components}

The calling sequence of this routine is
\begin{verbatim}
     laguerre(a,n,z,lagp,ierr)
\end{verbatim}
  where the input data are:  $a$, $n$ and $z$ (arguments of the Laguerre polynomial). 
  The outputs of the function are error flag $ierr$ and the value of the Laguerre polynomial value $lagp$.
  The possible values of the error flag are: $ierr=0$, successful 
computation; 
$ierr =1$, computation failed due to overflow/underflow; 
$ierr=2$, arguments out of range.

\section{Testing the algorithms}

The performance of the asymptotic expansions for the Laguerre polynomials has been tested 
by considering the relation given in Eq.(18.9.13) of  \cite{koorn:2010:OPS}
written in the form

\begin{equation}\label{eq:comp1}
\epsilon_1=\left|\Frac{L_{n-1}^{(\alpha+1)}(z)+L_n^{(\alpha)}(z)}{L_n^{(\alpha+1)}(z)}-1\right|.
\end{equation}

This chek fails close to the zeros of $L_n^{(\alpha+1)}$; in this case, we can consider the
alternative test

\begin{equation}\label{eq:comp2}
\epsilon_2=\left|\Frac{L_{n}^{(\alpha+1)}(z)-L_{n-1}^{(\alpha+1)}(z)}{L_n^{(\alpha)}(z)}-1\right|.
\end{equation}

Notice that, because the zeros of $L_n^{(\alpha)}$ and  $L_n^{(\alpha+1)}$
are interlaced, both tests will not fail simultaneously. We can therefore take

\begin{equation}\label{eq:comp}
 \epsilon=\mbox{min}(\epsilon_1,\epsilon_2).
\end{equation}
 
\begin{figure}
\begin{center}
\epsfxsize=15cm \epsfbox{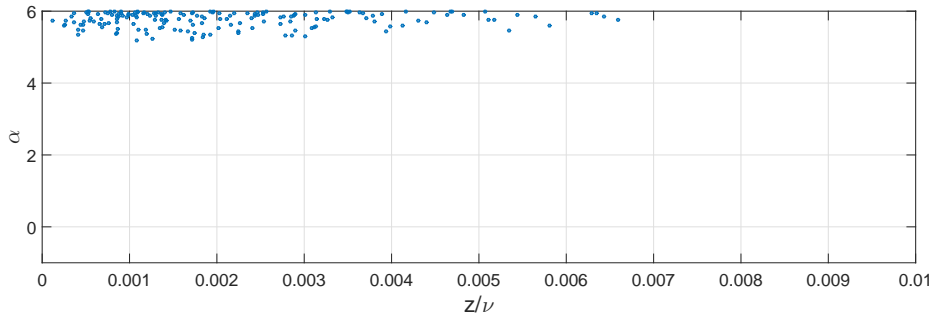}
\caption{Test of the performance of the simple Bessel-type expansion of section  (\ref{sec:Lagsim}). The points where the value of $\epsilon$ in Eq. (\ref{eq:comp}) is greater than
$5\,10^{-12}$ are plotted. 
\label{Fig2}}
\end{center}
\end{figure}

\begin{figure}
\begin{center}
\epsfxsize=15cm \epsfbox{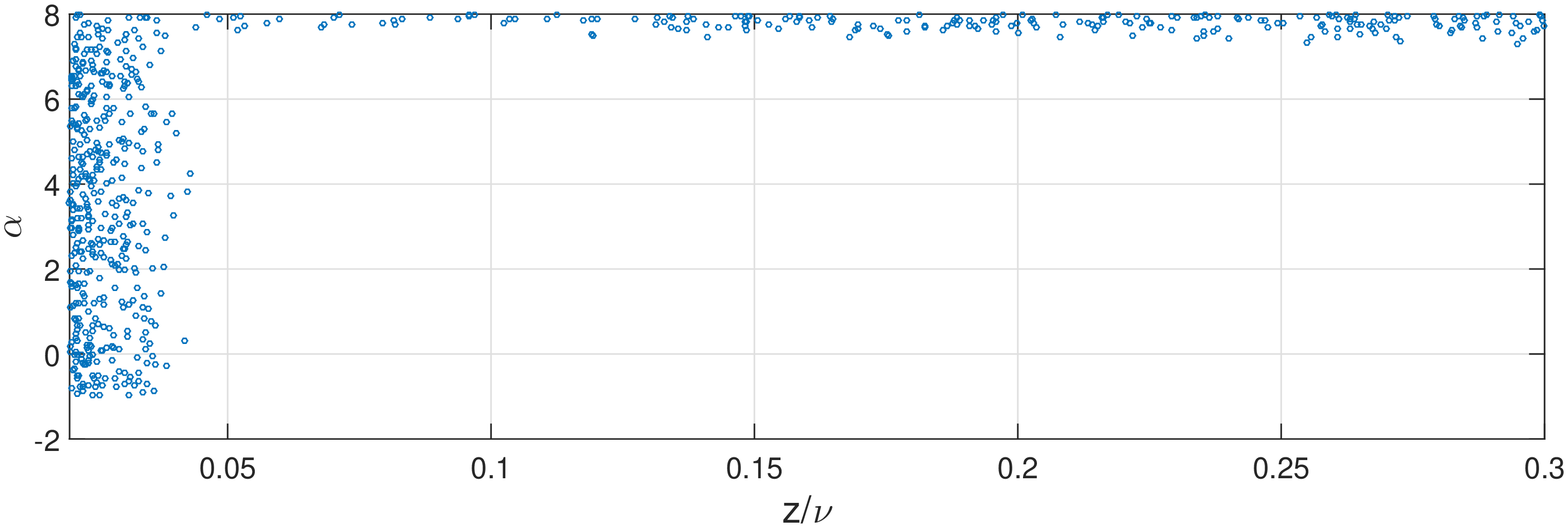}
\caption{Test of the performance of  the Bessel-type expansion of section (\ref{sec:LagBesexp}). The points where the value of $\epsilon$ in Eq. (\ref{eq:comp})  is greater than
$5\,10^{-12}$ are plotted.
\label{Fig3} }
\end{center}
\end{figure}

\begin{figure}
\begin{center}
\epsfxsize=15cm \epsfbox{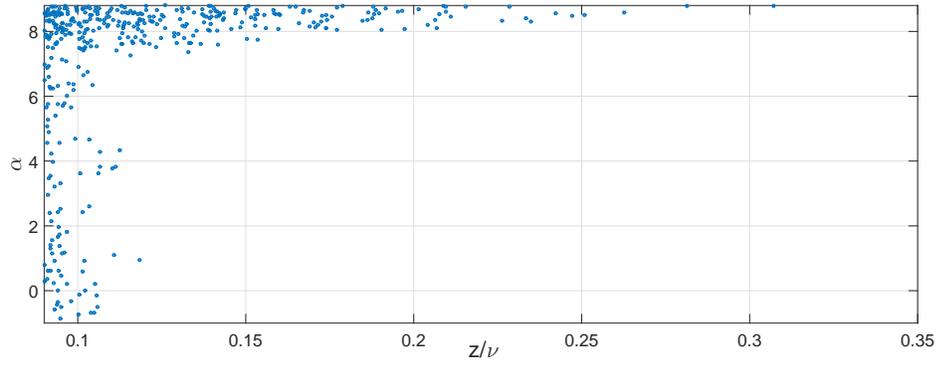}
\caption{Test of the performance of the Airy-type asymptotic expansion. The points where value of $\epsilon$ in Eq. (\ref{eq:comp})  is greater than
$5\,10^{-12}$ are plotted.
\label{Fig4} }
\end{center}
\end{figure}

A test for the accuracy obtained using the Bessel-type and Airy-type expansions for $n$ large is shown in Figures \ref{Fig2},  \ref{Fig3}
and \ref{Fig4}, respectively.
In these figures, the points where the value of $\epsilon$ in Eq. (\ref{eq:comp}) is greater than
$5\,10^{-12}$ are shown. Parameter values have been randomly generated in the oscillatory  region
of the Laguerre polynomials with $ (\alpha, \,n) \in (-1,\,10) \times (200,\,10000)$.

As can be seen, the Bessel-type expansions (Figures \ref{Fig2}, \ref{Fig3}) are valid for small values of the variable $z/\nu$; in particular, the simple
Bessel-type expansion works well when the argument of the Laguerre
polynomials is very close to the origin. On the other hand, the domain of applicability of the Airy-type expansion extends well beyond
the transition between the oscillatory and the monotonic regions of the functions, as Figure \ref{Fig4} shows for the oscillatory region. 
The validity of the three asymptotic
expansions is, in all cases, limited by the value of the parameter $\alpha$, as commented in \S(\ref{sec:LagBesexpla}). In the Fortran 90
module, we restrict the values of this parameter to the interval $(-1,\,5)$ in order to avoid the use of the recursion relation (\ref{ttrr})
for large $n$.

We have also tested the efficiency of using the asymptotic expansions in their region of applicability in comparison with
the use of the three-term recurrence relation given in Eq. (\ref{ttrr}) for computing the functions. 
Our tests show that for $n>200$ it is more efficient to use the asymptotic expansions than the recurrence relation.
As an example, Table 1 shows few of the test values for particular choices of the parameters.

\begin{table}
\label{table1}
$$
\begin{array}{cccc}
n  & z/\nu  & \mbox{CPU time AE} (s) &    \mbox{CPU time TTRR} (s)    \\
  \hline
 125 & 0.001 &  0.12   & 0.046     \\
       &   0.15       &  0.047    & 0.046     \\
      &     0.7     &  0.078    & 0.047       \\
200 &    0.001    &   0.047   &   0.078    \\
      &    0.15    &    0.047  & 0.078        \\
      &    0.7   &   0.078   &  0.078      \\
500 &    0.001    &   0.047   &   0.2    \\
      &    0.15    &   0.047   & 0.2       \\
      &    0.7   &    0.078  & 0.2       \\
1000 &    0.001    &  0.031    &   0.39    \\
      &    0.15    &   0.047   &  0.39      \\
10000 &    0.001    &  0.031    &   3.82    \\
\hline
\end{array}
$$
{\footnotesize {\bf Table 1}. Comparison of the efficiency of the asymptotic expansions (AE) vs the
use of the three-term recurrence relation (TTRR) given in Eq. (\ref{ttrr}) for computing the Laguerre polynomials
$L^{(1.5)}_n(z)$. Values of CPU times shown correspond to 20000 function evaluations. The parameter $\nu$ is
given in Eq. (\ref{eq:nu}).   }
\end{table}

\section{Test run description}

The Fortran 90 test program {\bf testlag.f90} includes the computation
of 25 function values and their comparison with the corresponding
pre-computed results. Also, the relation given in (\ref{eq:comp1})
is tested for several values of the parameters $(z,\,\alpha,\,n)$.

\section{Acknowledgements}

A.G. acknowledges the Fulbright/MEC Program for support during her stay at SDSU.
J.S.  acknowledges the Salvador de Madariaga Program for support during his stay at SDSU.
The authors acknowledge financial support from 
{\emph{Ministerio de Ciencia e Innovaci\'on}}, projects MTM2012-34787, MTM2015-67142-P. NMT thanks CWI, Amsterdam, for scientific support.




\section*{References}

\bibliographystyle{elsarticle-num}
\bibliography{laguerrebib}



\end{document}